\documentclass[reqno]{amsart}

\usepackage[all]{xy}

\numberwithin{equation}{section}

\theoremstyle{plain}
\newtheorem{lemma}{Lemma}[section]
\newtheorem{teo}[lemma]{Theorem}

\newtheorem{propo}[lemma]{Proposition}

\newtheorem{con}[lemma]{Conjecture}

\newtheorem{claim}{Claim}

\theoremstyle{definition}
\newtheorem{defi}[lemma]{Definition}
\newtheorem{defs}[lemma]{Definitions}
\newtheorem{example}[lemma]{Example}

\theoremstyle{remark}
\newtheorem{remark}[lemma]{Remark}
\newtheorem{notation}[lemma]{Notation}

\newcommand{\pic}{{\rm Pic\thinspace}}

\newcommand{\p}{\mathbb{P}}

\newcommand{\oc}{{\mathcal O}}

\newcommand{\ls}{{\mathcal L}}

\newcommand{\ci}{{\mathcal I}}

\newcommand{\cre}{{\rm Cr\thinspace}}

\newcommand{\rt}{\longrightarrow}
\newcommand{\rmap}{\dashrightarrow}

\begin{document}

\title{Elementary $(-1)$-curves of $\p^3$}
\author{Antonio Laface}
\address{
Antonio Laface \newline
Department of Mathematics and Statistics \newline
Queen's University \newline
K7L 3N6 Kingston, Ontario (Canada)
}
\email{alaface@mast.queensu.ca}

\author{Luca Ugaglia}
\address{
Luca Ugaglia \newline
Via Petrarca 24, \newline 14100 Asti, Italy }
\email{luca.ugaglia@gmail.com}
\keywords{Linear systems, fat points} \subjclass{14C20}
\begin{abstract}
In this note we deal with rational curves in $\p^3$ which are images 
of a line by means of a finite sequence of cubo-cubic Cremona transformations.
We prove that these curves can always be obtained applying to the line a sequence 
of such transformations increasing at each step the degree of the curve. 
As a corollary we get a result about curves that can give speciality 
for linear systems of $\p^3$.
\end{abstract}
\maketitle

\section*{Introduction}
Let $C$ be a rational integral curve contained in the blowing-up $X_n
$ of $\p^n$ along a set of points in general position. The notion of $
(-1)$-curve is defined only in case $n=2$ by asking the normal bundle
of $C$ to be $\oc_{\p^1}(-1)$. In order to find a good generalization
of $(-1)$-curves to the blowing-up of higher dimensional projective
spaces one is naturally lead to consider those curves whose normal
bundle is $\oc_{\p^1}(-1)^{n-1}$. It turns out that this definition
is not good enough since it can happen that $C$ is not unique inside
its rational equivalence class. In this paper we propose a definition
of a  generalized $(-1)$-curve of $\p^3$ which is based on the
observation that there exists a well defined action on $A^2
(X_n)$ of the group generated by monomial Cremona
transformations of bidegree $(3,3)$.
As in the case of $\p^2$ we define a $(-1)$-curve to be any element
which is in the orbit of the class of line through two simple points. \\
The construction of these curves (that we call elementary $(-1)$-
curves) has an application to the theory of linear systems.  The
system $\ls$ of surfaces of degree $d$ passing through $r$ general
points of $\p^3$ with multiplicities $m_1,\ldots, m_r$
is defined to be special if the conditions imposed by the multiple
points are dependent.
In~\cite{lu} two different types of special systems are constructed.
The first type is a
generalization of the $(-1)$-special systems of $\p^2$ defined in~
\cite{cm,har,hi,mir}.
The speciality of such a system $\ls$ depends on the fact that
there exists an elementary $(-1)$-curve $C$ such that $\tilde{\ls}C
\leq -2$, where $\tilde{\ls}$ is the strict transform of $\ls$.
In this paper we prove that given a linear system $\ls$ such that its
degree
cannot be decreased by means of one of the Cremona transformations described before,  then the only elementary
$(-1)$-curves that can give speciality are lines through two points.\\
The paper is organized as follows: in Section $1$ we provide some
preliminary material
and in Section $2$ we give the definition of elementary $(-1)$-curve
and we construct an example in which the definition doesn't  make
sense due to the lack of generality of the points. Section $3$
contains the numerical results needed for the proof of the two main
theorems of Section $4$ and in the final section we give some examples and a conjecture about elementary $(-1)$-curves.

\section{Preliminaries}

We start by fixing some definitions and notations.
\begin{defs}
Given a collection of points $q_1,\ldots,q_r\in\p^3$ and given a 
curve $C\subset\p^3$ (not necessarily irreducible) of degree $\delta$, passing through 
$q_i$ with multiplicity $\mu_i$, for $i=1,\ldots,r$, we define the {\em type of $C$} to be 
the sequence $(\delta;\mu_1,\ldots,\mu_r)$.
The set of all curves of type $(\delta;\mu_1,\ldots,\mu_r)$ through points in general position is denoted by 
$\ell(\delta;\mu_1,\ldots,\mu_r)$.

With $\ls(d;m_1,\ldots,m_r)$ we denote the linear system of 
surfaces of degree $d$ passing through the $q_i$'s with multiplicities $m_i$'s. 

Let us fix $r>0$ and let us consider the space of the configurations of $r$ points in $\p^3$, i.e. the quotient of $(\p^3)^r\setminus\Delta$ by the symmetric group $S_r$. For each choice of $d\geq 1$ and $m_i\geq 0$ (for $i=1,\ldots,r$) we denote by $U(d;m_1,\ldots,m_r)$ the set of configurations for which the dimension of the system $\ls(d;m_1,\ldots,m_r)$ is not the biggest possible. Let $U_r$ be the union of the closed subsets $U(d;m_1,\ldots,m_r)$, for each choice of $d$ and $m_i$. We say that the points $q_1,\ldots,q_r$ are {\em in general position} if their configuration does not lie in $U_r$.

From now on $\ls(d;m_1,\ldots,m_r)$ will denote a linear system with points in general position.

If the points $q_i$'s lie on a smooth quadric, we denote with $\ell_0(a,b;\mu_1,\ldots,\mu_r)$ 
the linear system of curves $\oc(a,b)$ through the multiple points.

\end{defs}
Let us denote by $X$ the blowing up of $\p^3$ along the points $q_i$.
Given $\ls=\ls(d;\allowbreak m_1,\allowbreak \ldots,\allowbreak m_r)$ and $C\in\ell(\delta;
\mu_1,\ldots,\mu_r)$,
by abuse of notation we denote by $\ls C$ the intersection product of the strict transforms
$\tilde{\ls}$ and $\tilde{C}$ on $X$, i.e. $\ls C=d\delta-\sum_{i=1}^r m_i\mu_i$.

Given four non-planar points $q_1,\ldots,q_4$ in $\p^3$, let us consider a cubo-cubic Cremona 
transformation based on these points and expressed in coordinates in the following way: 
$(x_0:\ldots :x_3)\rmap (x_0^{-1}:\ldots :x_3^{-1})$. 
\begin{notation}
Throughout the paper, by abuse of notation we will say that a curve $C$ does not intersect the $1$-dimensional indeterminacy locus of a cubo-cubic Cremona transformation if either $C$ does not intersect any of the six lines through the base points or if the intersection is along these points and transversal.
\end{notation}
The action of a cubo-cubic Cremona 
transformation on linear systems and sets of curves which do not intersect the $1$-dimensional indeterminacy locus out from the base points is given by (see~\cite{lu}):
\begin{eqnarray}\label{cre-a1}
\ls(d;\allowbreak m_1,\allowbreak \ldots,m_r) & \mapsto & \ls(d+k;m_1+k,\ldots, m_4+k,m_5,
\ldots,m_r),\\ \label{cre-a2}
\ell(\delta;\allowbreak \mu_1,\allowbreak
\ldots,\mu_r) & \mapsto & \ell(\delta+2\gamma;\mu_1+\gamma,\ldots,
\mu_4+\gamma,\mu_5,\ldots,\mu_r),
\end{eqnarray}
where $k=2d-\sum_{i=1}^4 m_i$ and $\gamma=\delta-\sum_{i=1}^4 \mu_i$. 

Let us recall here some definitions from~\cite{lu}:

\begin{defs}
An {\em elementary transformation} is a finite composition of cubo-cubic Cremona transformations such that the total set of points involved is in general position.

A system $\ls(d;m_1,\ldots,m_r)$ is said to be in {\em standard form} if 
$m_1\geq\ldots\geq m_r$ and $d\geq\sum_{i=1}^4m_i$
(i.e. it is not possible to decrease its degree by means of an elementary transformation).
\end{defs}

\section{Elementary $(-1)$-curves}\label{-1c}

Let us consider an elementary transformation $\varphi=\sigma_s\circ\cdots\circ\sigma_1$ based on $q_1,\ldots,q_r$, points in general position of $\p^3$, and let us put $C_0=\langle q_1,q_2\rangle$. Throughout these notes, by induction we will denote $C_i$ the image of $C_{i-1}$ by $\sigma_i$, for $i=1,\ldots,s$. Let us give the following:

\begin{defi}
An {\em elementary $(-1)$-curve} is the image of the line $C_0$ via an elementary transformation $\varphi=\sigma_s\circ\cdots\circ\sigma_1$ such that for each $i=1,\ldots,s$ the curve $C_{i-1}$ does not intersect the $1$-dimensional indeterminacy locus of $\sigma_i$.
\end{defi}

We remark that if $C_{i-1}$ intersects the $1$-dimensional indeterminacy locus of $\sigma_i$ out from the $4$ base points, then  the action of $\sigma_i$ on the curve is different from the one given in equation~\eqref{cre-a2}, as we can see in the following example (see~\cite{lu} for a more detailed description of the action on the curves in this particular case).

\begin{example}
Let us take the line $C_0$ through the points $q_1$ and $q_2$ and let us consider a cubo-cubic transformation $\sigma$ based on the four points $q_3,\ldots,q_6$ such that $C_0$ intersects the line $\langle q_3,q_4\rangle$ in a point different from $q_i$ for $i=1,\ldots,4$. 
Let us denote by $\pi$ the plane $\langle C_0,q_3,q_4\rangle$ and by $q$ the intersection of $\pi$ with the line $\langle q_5,q_6\rangle$. If we restrict $\sigma$ to $\pi$ we get a quadratic Cremona transformation based on $q_3,q_4$ and $q$. The line $C_0$ can be seen as an element of the linear system $\ls_2(1;1^2)$, which is transformed into an element of $\ls_2(2;1^5)$.
This means that in this particular case the image of a line is not the rational cubic through $6$ fixed points (as we expect from formula~\eqref{cre-a2}), but it is a conic through $5$ planar fixed points.\\
Let us remark that in this case the points $\{q_1,\ldots,q_6\}$ are not in general position, since the first four lie on a plane.
\end{example}

In this section we are going to prove that in fact, if the points we are choosing are in general position then this phenomenon cannot happen. This is equivalent to say that if we apply an elementary transformation based on points in general position to a line through two points we do get an elementary $(-1)$-curve. In particular we can always apply equation~\eqref{cre-a2} and we deduce that every elementary $(-1)$-curve has odd degree. 

In order to prove this result we are going to use the following strategy.
We specialize the points $q_i$ on a quartic elliptic curve $E\subset\p^3$ and we consider the quadric $Q$ containing both $E$ and the line $C_0=\langle q_1,q_2\rangle$. Under these assumptions, we will prove that if we transform the line $C_0$ under an elementary transformation $\varphi=\sigma_s\circ\cdots\circ\sigma_1$, based on the $q_i$'s, then at each step $C_{i-1}$ does not intersect the $1$-dimensional indeterminacy locus of $\sigma_i$. Since this is true when we specialize the points $q_i$, it must hold also when the points are in general position in $\p^3$. 

We start by proving the following useful lemma about the restriction of a cubo-cubic Cremona transformation to a smooth quadric.

\begin{lemma}\label{crequa}
Let $\sigma$ be a cubo-cubic Cremona transformation based on the points $q_1,\ldots,q_4$, let $Q$ be a smooth quadric containing these four points and let $Q'=\sigma(Q)$. 
Therefore $\sigma$ induces an isometry of $\pic(\tilde{Q})$ and $\pic(\tilde{Q}')$ that can be represented by the following:
\begin{equation}\label{sigmaq}
\ell_0(a,b;\mu_1,\ldots \mu_4)\mapsto \ell_0(a+k,b+k;\mu_1+k,\ldots, \mu_4+k),
\end{equation}
where $k=a+b-\sum_{j=1}^4\mu_j$.
\end{lemma}

\begin{proof}
Let us consider the resolution of the 
indeterminacy of $\sigma$:
\[
\xymatrix{
\hspace{7mm} Y\ar[r]^-{\sigma_Y}\ar@<2.5ex>[d]_{p} & Y'\ar@<-5ex>[d]^{p'}
\hspace{12mm} \\
\tilde{Q}\subset X\ar@{-->}[r]^-{\sigma_X}\ar@<2.5ex>[d]_{\pi'} & X' 
\supset\tilde{Q}'\ar@<-5ex>[d]^{\pi'} \\
Q\subset \p^3\ar@{-->}[r]^-{\sigma} & \p^3 \supset Q', \\
}
\]
where $\pi$ is the blowing up of $q_1,\ldots,q_4\in Q_1$ (resp. $\pi'$ is the blowing up of $Q'$ along $q_1',\ldots,q_4'$),
and $p$ (resp. $p'$) is the blowing up of the strict transforms of the six lines through the preceding points. 
Since $Q$ does not contain any of these lines, the restriction of 
$\sigma_X:=p'\circ\sigma_Y\circ p^{-1}$ to $\tilde{Q}$ is a morphism and it
can be seen as a base change in 
$\pic(\tilde{Q})$. In order to see this, let $\{h_1,h_2,e_1,\ldots,e_4\}$ be a base of 
$\pic(\tilde{Q})$, where $h_1$ and $h_2$ are the pull-back of the rulings of $Q$ while
the $e_j$'s are the exceptional divisors corresponding to $q_1,\ldots q_4\in Q$.
We know by~\eqref{cre-a2} that a line is transformed by $\sigma$ into a rational 
normal cubic through $q'_1,\ldots,q'_4$. This implies that $h_1$ and $h_2$ are 
transformed into $2h_1+h_2-\sum_{j=1}^4 e_j$ and $h_1+2h_2-\sum_{j=1}^4 e_j$ respectively. 
Let $f_1,\ldots,f_4$ 
be the strict transforms of the four conics through three of the four points, i.e. 
$f_j \in |h_1+h_2-\sum_{k\neq j}e_k |$. Since each plane through three of $q_1,\ldots,q_4$ 
is collapsed by $\sigma$ to a point, each of these conics is contracted by 
${\sigma}_{|Q}$ and hence the base change exchanges $e_j$ with $f_j$. 
Therefore a curve $ah_1+bh_2-\sum_{j=1}^4e_j$ is transformed into $a(2h_1+h_2-\sum_{j=1}^4 e_j)+
b(h_1+2h_2-\sum_{j=1}^4 e_j)-\sum_{j=1}^4 \mu_j(h_1+h_2-\sum_{k\neq j}e_k)$, and we get the thesis.
\end{proof}

In what follows we denote by $E\subset\p^3$ a smooth elliptic quartic
and by $Q$ a smooth quadric containing $E$. With $D_{a,b}$ we will
denote the divisor $\oc_Q(a,b)_{|E}$ of $\pic(E)$.

\begin{lemma}\label{lem1}
There exist points $q_1,\ldots,q_r\in E$ and a quadric $Q$
such that $D_{1,0}-q_1-q_2$ is effective and $D_{a,b}-\sum m_iq_i\in
\pic^0(E)$ is effective only if $m_3=\ldots =m_r=0$.
\end{lemma}
\begin{proof}
Let $Q$ be the quadric of the pencil containing $\langle q_1,q_2
\rangle$. By moving $q_1$ if necessary,
we may assume that $Q$ is smooth.
Let us fix $a,b,m_1,\ldots,m_r$ such that $2(a+b)-\sum_{i=1}^rm_i=0$.
After reordering the points we may always suppose that $m_3\geq\ldots
\geq m_r$ so that if $m_3>0$, the map
\[
z\mapsto D_{a,b}-\sum_{i\neq 3}m_iq_i-m_3z.
\]
is defined. This map is an isogeny of $E$ into $\pic^0(E)$ and in
particular $0$ has a finite number of inverse images.
Therefore given $a,b,m_1,\ldots,m_r$ as before, the zero degree divisor
$D_{a,b}-\sum_{i\neq 3}m_iq_i-m_3z$ is not effective as $z$ varies in
a Zariski open set $U(a,b,m_1,\ldots,m_r)\subset E$.
Let $U=\cap U(a,b,m_1,\ldots,m_r)$ as the $r+2$-tuple varies over all
the admissible values, then by choosing $q_3\in U$ we get the thesis.
\end{proof}

\begin{lemma}\label{lem2}
Let $q_1,\ldots,q_r$ be as in Lemma~\ref{lem1} and let $\tilde{Q}
\stackrel{\pi}\rt Q$ be the blowing up of $Q$ along $q_1,\ldots,q_r$.
Then $C\in |\pi^*\oc_Q(1,0)-e_1-e_2|$ is the only rational curve in $
\tilde{Q}$ with  $C^2\leq -2$.
\end{lemma}
\begin{proof}
Let $C\in |\pi^*\oc_Q(a,b)-\sum_{i=1}^rm_ie_i|$ be a rational curve
with $C^2\leq -2$.
By construction the anticanonical divisor $-K_{\tilde{Q}}$ is
effective and linearly equivalent to $E$. This means that $C^2=-CK_
{\tilde{Q}}-2\geq -2$ and therefore $C^2=-2$. Moreover $\oc_E(C)\sim
\oc_E$ and by Lemma~\ref{lem1} one has that  $m_3=\ldots=m_r=0$. An
easy computation shows that $a+b=1$.
\end{proof}

Let us consider now an elementary transformation $\varphi=\sigma_s
\circ\cdots\circ\sigma_1$ based on $q_1,\ldots,q_r$ as in Lemma~\ref
{lem1} and let $C_i$ be as before. We suppose that
$C_{i-1}$ is not in the indeterminacy locus of $\sigma_i$ for $i=1,
\ldots,s$ and we denote $Q_0=Q$ and $Q_i=\sigma_i(Q_{i-1})$. The we
have the following:

\begin{lemma}\label{lem3}
$C_i\subset Q_i$ is the only rational curve of self intersection $
\leq -2$.
\end{lemma}
\begin{proof}
In the preceding lemma we have seen that the line $C_0$ through $q_1$
and $q_2$ is the only rational curve of self intersection $\leq -2$
on $Q_0$. By Proposition~\ref{crequa} $\sigma_i:\pic(\tilde{Q}_{i-1})\rt\pic(\tilde{Q}_i)$ is an isometry of lattices which preserves the effective cone, and hence we get the thesis.
\end{proof}

\begin{propo}
Let $\varphi=\sigma_n\circ\cdots\circ\sigma_1$ be an elementary transformation based on $q_1,\ldots,q_r$ as
before and such that $C_0=\langle q_1,q_2\rangle$ is not contained in the indeterminacy locus of $\varphi$.
Then $\varphi(C_0)$ is an elementary $(-1)$-curve.
\end{propo}
\begin{proof}
As before, let $C_i=\sigma_i(C_{i-1})$ and let us suppose that $C_i$ is an elementary $(-1)$-curve. Then 
$C_i\subset Q_i$ and $\tilde{C}_i^2=-2$ in $\pic(\tilde{Q}_i)$. We need only to check that $\sigma_{i+1}(C_i)$ 
is an elementary $(-1)$-curve. If this is not the case then $C_i$ intersects one of the six fundamental lines of 
$\sigma_{i+1}$ outside the four fundamental points. Let us call $L$ this line through two of the $q_i$'s, then $L\in Q_i$
(since it intersects $Q_i$ at least in three points) and moreover $\tilde{L}^2\leq -2$ in $\pic(\tilde{Q}_i)$. By lemma~\ref{lem3}
we have that $L=C_i$ which is a contradiction.
\end{proof}

The generality assumption made on the position of the points is important and can not be avoided. 
In order to see this one can consider the following example:
\begin{example}
Consider an elliptic quartic $E\subset\p^3$. 
\begin{claim}\label{claim}
There exist $p_1,\ldots,p_8\in E$ such that, for any $Q\in|\oc_{\p^3}(2)\otimes\ci_E|$, the divisor 
${K_Q}_{|E} + p_1+\cdots+p_8$ is a $2$-torsion point of $\pic^0(E)$. 
\end{claim}
From now on $Q$ will be a quadric of the pencil through $E$. 
Let $\pi: X\rt\p^3$ be the blowing-up of $\p^3$ along eight points with the claimed property and let 
$E_i=\pi^{-1}(p_i)$ and $e_i=E_i\cap\tilde{Q}$ where $\tilde{Q}$ is the strict transform of $Q$.
Observe that $|-K_{\tilde{Q}}|$ contains only one curve which is isomorphic to $E$. 
By abuse of notation let us call this curve $E$. Consider the exact sequence:
\[
\xymatrix{
0\ar[r] & \oc_{\tilde{Q}}(E)\ar[r] & \oc_{\tilde{Q}}(2E)\ar[r] & \oc_{E}(2E) \ar[r] & 0,
}
\]
then $\dim |2E| = 1$ since
$\oc_E(2E)=-2(\pi^*{K_Q}+e_1+\cdots+e_8)_{|E}\sim \oc_E$ by hypothesis. \\
Consider now the linear system $|2\tilde{Q}|$ and observe that its dimension can be evaluated 
by means of the exact sequence
\[
\xymatrix{
0\ar[r] & \oc_{X}(\tilde{Q})\ar[r] & \oc_{X}(2\tilde{Q})\ar[r] & \oc_{\tilde{Q}}(2\tilde{Q}) \ar[r] & 0,
}
\]
and the fact that $\oc_{\tilde{Q}}(2\tilde{Q}) = 2E$. From the preceding sequence we see that
$\dim |2\tilde{Q}| = 3$ which in particular means that there exists an irreducible quartic surface $S$
which is singular along $p_1,\ldots,p_8$. \\
Consider now an elementary $(-1)$-curve $C$ of type $(\delta;\mu_1,\ldots,\mu_8)$ through the
$p_i$'s. Observe that $C\tilde{Q}=0$, so that given $p\in C$ the element of $|\tilde{Q}|$ which contains $p$
contains also $C$. By abuse of notation let us call this element still $\tilde{Q}$. In the same way let $S$ be 
an element of $|2\tilde{Q}|$ which contains $C$. From the preceding discussion we have that
$S_{|\tilde{Q}} = 2E$ which implies that 
\[
2E = C + C'
\]
where $C'$ is an effective curve of $\tilde{Q}$. This fact implies that the possible types for $C$ are finite, since its degree must be smaller than $8$ (which is $\deg 2E$).
Moreover, since $|2E|$ is an elliptic pencil and $C$ is a $(-2)$-curve on $\tilde{Q}$, then 
$C'$ is union of $(-2)$-curves. 
\end{example}

\begin{proof}[Proof of Claim~\ref{claim}]
Given $D_t = {K_{Q_t}}_{|E}$, where $Q_t$ is a quadric of the pencil containing $E$ and parametrized by $t$, the map $t\mapsto [D_t+p_1+\cdots+p_8]$
of $\p^1$ into $\pic^0(E)$ is necessarily constant. We conclude by observing that for a fixed value of $t$, say $t=0$, the map $z\mapsto [D_0+p_1+\cdots+p_7+z]$
is an isomorphism of $E$ and $\pic^0(E)$.
\end{proof}

\section{Some numerical results}

In this section we give some characterizations of elementary
$(-1)$-curves that we are going to use later in order to prove our main result.
\begin{defi}
Given two curves $C\in\ell(\delta;\mu_1,\ldots,\mu_r)$ and $C'\in\ell(\delta';\mu_1',\ldots,
\mu_r')$, we define the bilinear form $\langle\ ,\ \rangle$ as 
\[
\langle C,C'\rangle = \delta\delta'-2\sum_{i=1}^r\mu_i\mu_i',
\]
which is invariant under elementary transformations.
\end{defi}
Since an element $K\in\ell(4;1^r)$ is also invariant, we get that  
$\langle C,C\rangle$ and $\langle C,K\rangle$
are preserved by elementary transformations.
In particular if $C$ is an elementary $(-1)$-curve, then 
\begin{equation}\label{-1}
\begin{array}{llc}
\langle C,C\rangle & = & -3\\ 
\langle C,K\rangle & = & 0, 
\end{array}
\end{equation}
since these are the values for the line $\ell(1;1^2)$. 
\begin{defi}
We say that a curve $C\in\ell(\delta;\mu_1,\ldots,\mu_r)$ is a 
{\em numerically elementary $(-1)$-curve} if it is irreducible, reduced and satisfies~\eqref{-1}. 
\end{defi}

We end the section with another numerical condition satisfied by the elementary $(-1)$-curves, but in order to do this we need the following:

\begin{lemma}\label{qua}
Let $C\in\ell(\delta;\mu_1,\ldots,\mu_r)$ be an elementary $(-1)$-curve. Then there exists a curve $C'$ having the same type of $C$,
contained in a smooth quadric and such that $C'\in\ell_0(\alpha+1,\alpha;\mu_1,\ldots,\mu_r)$ 
(where $\delta=2\alpha+1$). 
\end{lemma}
\begin{proof}
Let $\varphi=\sigma_s\circ\cdots\circ\sigma_1$ be the elementary transformation which 
sends the line $C_0$ to $C$. Let us take a quadric $Q\supset C_0$ and consider 
the cubo-cubic Cremona transformation $\sigma_1'$ based on four general 
points of $Q$. By~\eqref{cre-a1}, $Q_1:=\sigma_1'(Q)$ is a smooth quadric and 
$\sigma_1'(C_0)$ is a rational curve contained in $Q_1$. 
Proceeding with this idea, we choose the fundamental points of $\sigma_2'$ on $Q_1$ and 
so on. In this way, after $n$ steps we have constructed an elementary transformation
$\varphi'=\sigma_s'\circ\cdots\circ\sigma_1'$ such that $Q':=\varphi'(Q)$ is a 
smooth quadric and $C':=\varphi'(C_0)$ is an elementary $(-1)$-curve of the same type of $C$ and lying on $Q'$. Therefore we can write $C'\in\ell_0(a,b;\mu_1,\ldots \mu_r)$, whith $a+b=\delta$, and by Lemma~\ref{crequa} we have that the difference $a-b$ is preserved by each $\sigma_i'$. Since $C'$ is the image
of the line $C_0\in\ell_0(1,0;1^2)$, we have that $a-b=1$, which proves the last part of the lemma. 
\end{proof}

\begin{propo}\label{d>2m1}
If $C\in\ell(\delta;\mu_1,\ldots,\mu_r)$ is an elementary $(-1)$-curve through $r$ points in 
general position then either $C\in\ell(1;1^2)$ or $\delta > 2\mu_1$.
\end{propo}
\begin{proof}
By the Lemma~\ref{qua} we know that $C$ has the same type of 
$C'$ contained in a smooth quadric $Q$. Moreover 
$C'\in\ell_0(\alpha+1,\alpha;\mu_1,\ldots,\mu_r)$ 
and it is irreducible since it is the image of a line by a finite sequence of 
cubo-cubic Cremona transformations. Let us suppose that 
$\delta\leq 2\mu_1-1$ (since $\delta$ is odd, it cannot be $2\mu_1$).
Then $C'(h_1-e_1)=\alpha-\mu_1 < 0$ (we make the intersection product on the blow up
of the quadric), which means that the line $h_1-e_1$ is a 
component of $C'$. But since $C'$ is irreducible, $C'= h_1-e_1$.
\end{proof}

\section{A theorem on elementary $(-1)$-curves}
The aim of this section is to prove that given an elementary $(-1)$-curve $C$, it is 
always possible to find a sequence of elementary transformations $\rho_1,\ldots,\rho_n$ 
such that, if we put $C_n = C$ and $C_{i-1} = \rho_{i}(C_{i})$ for $i=n,\ldots,1$, then the sequence $\{C_i\}$ satisfies:  
$C_0 = \ell(1;1^2)$ and $\deg C_{i} < \deg C_{i+1}$. In other words this means that there exists a sequence of elementary transformations which sends $C$ to a line through two points and such that the degree decreases at each step. \\
This fact is used later in order to prove that (Theorem~\ref{mt}), given a linear system 
$\ls$ of $\p^3$ in standard form, the only elementary $(-1)$-curves which can have negative 
intersection with $\ls$ are lines through two points. \\
In order to prove the main result we will use the numerical 
characterizations of elementary $(-1)$-curves. 
The following lemma guarantees that if we apply a cubo-cubic Cremona transformation to 
an elementary $(-1)$-curve (different from the line), then the type of the new curve 
can not contain negative numbers. 
\begin{lemma}
Let $C$ be an elementary $(-1)$-curve and let $C'=\sigma(C)$ be the image of $C$ 
by means of a Cremona transformation. Then $C'$ can have negative degree or a 
negative multiplicity if and only if $C\in\ell(1;1^2)$.
\end{lemma}
\begin{proof}
Let $C\in\ell(\delta;\mu_1,\ldots,\mu_r)$ and let us reorder the points in order to
have $\mu_1\geq\ldots\geq\mu_r$. We can reduce to the case in which 
$\sigma$ is the Cremona transformation based on the first four points $q_1,\ldots,q_4$.
We denote by $C'\in\ell(\delta';\mu_1',\ldots,\mu_s')$ the image of $C$ and we 
first suppose that $\mu_4' < 0$. By equation~\eqref{cre-a2} this is equivalent to say that 
$\delta - \sum_{i=1}^3\mu_i < 0$, and hence the curve $C$ intersects negatively
the plane $\langle q_1,q_2,q_3\rangle$. Since $C$ is irreducible, it must be 
contained in that plane. Moreover we deduce that $r=3$, since the points 
$q_i$ are general. 
By means of equations~\eqref{-1}, rearranging the terms we get that 
$\sum_{i < j}(\mu_i-\mu_j)^2 + 5\sum_{i=1}^3\mu_i^2 = 12$. But the only integer solution
of this equation is $\mu_1=\mu_2=1$ and $\mu_3=0$, corresponding to the line
through $q_1$ and $q_2$.\\
Let us suppose now that $\delta' < 0$ and $C$ is not a line through two points. 
By equation~\eqref{cre-a2} we have that $3\delta - 2\sum_{i=1}^4\mu_i < 0$ and by 
the preceding lemma $\delta > 2\mu_1$. Therefore we get that 
$\delta < \sum_{i=2}^4\mu_i$, which means that the curve $C$ intersects negatively 
the plane $\langle q_2,q_3,q_4\rangle$. Arguing as before we get a contradiction.
\end{proof}
\begin{remark}
We remark that when we take a line $L$ passing through two points $q_1$ and $q_2$, 
and we apply a cubo-cubic Cremona transformation based on $q_1,\ q_2$ and two other
points out from the line, numerically we get a curve in $\ell(-1;0^2,-1^2)$. 
Actually, this curve does 
not exist, since the line $L$ is contained in the indeterminacy locus of the 
transformation.
\end{remark}
We are now able to prove the following:
\begin{teo}\label{mainpropo}
Any elementary $(-1)$-curve can be obtained from $\ell(1;1^2)$ by a finite sequence 
of cubo-cubic Cremona transformations such that at each step the degree increases.
\end{teo}
\begin{proof}
Let us consider an elementary $(-1)$-curve $C\in\ell(\delta;\mu_1,\ldots,\mu_r)$, different 
from a line through $2$ points and let us reorder the $q_i$'s in order to have 
$\mu_1\geq\mu_2\geq\ldots\geq\mu_r$. 
We claim that under our hypothesis we have $\delta<\sum_{i=1}^4 \mu_i$.\\
Since $C$ is an elementary $(-1)$-curve, by~\eqref{-1} we have that
\[
\langle C,C\rangle + \mu_4\langle C,K\rangle < 0.
\]
This is equivalent to say that $\delta^2-2\sum_{i=1}^r\mu_i^2+\mu_4
(2\sum_{i=1}^r\mu_i-4\delta) < 0$,
and hence
\[
\delta(\delta-4\mu_4)-2(\sum_{i=1}^4\mu_i^2-\mu_4\sum_{i=1}^4\mu_i) < 2(\sum_{i=5}^r\mu_i^2 - 
\mu_4\sum_{i=5}^r\mu_i ) 
\leq 0,
\]
where the last inequality holds since $\mu_4\geq \mu_i$, for $i\geq 5$.
By Proposition~\ref{d>2m1} we know that $\delta > 2\mu_1$; let us suppose by contradiction 
that $\delta\geq \sum_{i=1}^4\mu_i$.
By substituting these two values in the preceding inequality we obtain:
\[
2\mu_1(\sum_{i=1}^4\mu_i-4\mu_4)-2(\sum_{i=1}^4\mu_i^2-\mu_4\sum_{i=1}^4\mu_i) < 0,
\]
which is equivalent to say 
\[
(\mu_1-\mu_2)(\mu_2-\mu_4)+(\mu_1-\mu_3)(\mu_3-\mu_4) < 0,
\]
a contradiction which proves the claim.\\
Now we can proceed by induction. Let us denote by $\sigma$ the cubo-cubic Cremona 
transformation based on the first four points, and let $C'\in\ell(\delta';
\mu_1',\ldots,\mu_r')$ be the image $\sigma(C)$. By the preceding lemma, $\delta'$ 
and $\mu_i'$ are non-negative. Moreover, since $\delta < \sum_{i=1}^4\mu_i$,
by equation~\eqref{cre-a2} we have that $\delta' < \delta$. 
If $C'$ is a line through two points we are done. Otherwise, reordering the 
multiplicities, by Proposition~\ref{d>2m1} we have $\delta' > 2\mu_1'$ and 
we can prove as before that $\delta' < \sum_{i=1}^4\mu_i'$. We can then apply
a new cubo-cubic Cremona transformation decreasing the degree, and so on.
\end{proof}
As a corollary of this result we finally prove the following theorem concerning
linear systems in standard form (see~\cite[Corollary 5.3]{lu}):
\begin{teo}\label{mt}
If $\ls=\ls(d;m_1,\ldots,m_r)$ is a non-empty linear system in standard form and $C$ 
is an elementary $(-1)$-curve such that $\ls C< 0$, then $C\in\ell(1;1^2)$.
\end{teo}
\begin{proof}
Let us consider an elementary $(-1)$-curve $C$ which is not in $\ell(1;1^2)$.
By Proposition~\ref{mainpropo} we can suppose that $C$ is obtained from $\ell(1;1^2)$
by a finite set of cubo-cubic Cremona transformations increasing the degree.
We claim that under these hypotesis we can write 
\begin{equation}\label{int}
\ls C=\beta_1(2d-\sum_{i=1}^4m_i^{(1)})+\cdots+
\beta_a(2d-\sum_{i=1}^4m_i^{(a)})+(d-m_h-m_k),
\end{equation}
where $\beta_j\geq 1$, $m_i^{(j)}$ are chosen between $m_1,\ldots,
m_r$ and $h,k\geq 5$.\\
We argue by induction on the number of Cremona transformations
necessary to obtain $C$ from the line $\ell(1;1^2)$. First of
all, after one transformation the image of the line is the
rational normal cubic $\ell(3;1^6)$, having intersection product
with $\ls$ equal to $3d-\sum_{i=1}^6 m_i=(2d-\sum_{i=1}^4
m_i)+(d-m_5-m_6)$. Now we assume that the formula is true for
$C\in\ell(\delta;\mu_1,\ldots,\mu_s)$ and we prove it for the
curve $C'=\cre(C)\in\ell(\delta';\mu_1',\ldots,\mu_s')$
obtained from $C$ performing one more Cremona transformation
increasing the degree. Reordering the multiplicities we can suppose that the transformation is
based on the first $4$ points. By formula~\eqref{cre-a2},
$\delta'=\delta + 2\gamma$ and $\mu_i'=\mu_i + \gamma$ for $i=1,\ldots,4$,
where $\gamma=\delta-\sum_{i=1}^4 \mu_i>0$, and $\mu_i'=\mu_i$ for
$i\geq 5$. Therefore $\ls C'-\ls C=\gamma(2d-\sum_{i=1}^4 m_i)$,
which proves the claim.\\
Since $\ls$ is in standard form, $2d\geq\sum m_i^{(j)}$ and
$d\geq m_h + m_k$ (otherwise $2d < m_1 + m_2 + m_h + m_k$), and
hence all the terms on the right side of equation~\ref{int} are
non-negative.
\end{proof}

The following example shows that Theorem~\ref{mt} is no longer true in higher dimension.

\begin{example}
As in the case of $\p^3$ we can consider the elementary transformation of $\p^4$ given 
by $(x_0:\ldots:x_4)\rmap(x_0^{-1}:\ldots:x_4^{-1})$, and say that a linear system is 
in standard form if its degree cannot be decreased by means of such a transformation.\\
Let us fix $q_1,\ldots,q_7\in\p^4$ in general position.
The system $\ls$ of hypersurfaces of degree $5$ with multiplicity $3$ in $q_i$,
for $i=1\ldots,7$, is not empty since its virtual dimension is positive. 
In~\cite{lu1} it is proved that by applying to $\ls$ the preceding transformation, 
based on five of the $q_i$'s, it remains unchanged. In particular this implies 
that $\ls$ is in standard form. Moreover the line through the remaining 
two points is transformed into the rational normal quartic $C$ through 
$q_1,\ldots,q_7$. This means that $C$ is an elementary $(-1)$-curve and $\ls C=-1$.
\end{example}

\section{A conjecture about elementary $(-1)$-curves}

In this section we will give some examples of sequences $(\delta;\mu_1,\ldots,\mu_r)$ 
satisfying equalities~\eqref{-1} but which do not correspond to elementary
$(-1)$-curves. 
These examples, together with Proposition~\ref{proj}, suggest
the following:

\begin{con}\label{conj1}
The only numerically elementary $(-1)$-curves satisfying $\delta < 2\mu_1$ 
are the lines through two points.
\end{con}

Let us consider a numerically elementary $(-1)$-curve $C$, different from the line. Then
by Conjecture~\ref{conj1}, $\delta > 2\mu_1$. Following the proof of Theorem~\ref{mainpropo} 
we know that there exists a cubo-cubic transformation $\sigma$ decreasing the degree of $C$ 
and the curve $C'=\sigma(C)$ has the same numerical properties of $C$ unless it is a line 
through two points.
So, after a finite number of transformations, we obtain a line through two points, 
which is equivalent to say that $C$ is an elementary $(-1)$-curve. Therefore 
Conjecture~\ref{conj1} is equivalent to the following:

\begin{con}\label{conj2}
Every numerically elementary $(-1)$-curve is an elementary $(-1)$-curve.
\end{con}

In order to give some evidence for Conjecture~\ref{conj1} we are going to prove 
that (assuming Harbourne-Hirschowitz Conjecture for linear systems on $\p^2$), 
there are no numerically elementary $(-1)$-curve different from
a line and satisfying $\delta = 2\mu_1 - 1$. 

\begin{propo}\label{proj}
Harbourne-Hirschowitz conjecture implies that if $\delta =2\mu_1-1 > 1$, then there
are no numerically elementary $(-1)$-curves of type $(\delta;\mu_1,\ldots,\mu_r)$. 
\end{propo}
\begin{proof}
It is enough to prove that a curve $C$ of type $(2\mu_1-1;\mu_1,\ldots,\mu_r)$ such that
$\langle C,C\rangle = -3$, $\langle C,K\rangle = 0$ and $\mu_1 > 1$, must be reducible 
or non-reduced.
Let us denote by $q_i\in C$ the points of multiplicity $\mu_i$ and consider the projection
$\pi$ of $C$ from $q_1$ to a generic plane $\Pi$. The image $\pi(C)$ is given 
by $\Gamma\cup\{q'_1,\ldots,q'_s\}$, where $\Gamma$ is a plane curve of degree $\mu_1 - 1 > 0$ 
and the $q'_i$ are points corresponding to lines through $q_1$. We can suppose that there 
are no such points, since otherwise $C$ is reducible. Then 
$\Gamma$ belongs to the linear system $\ls_2=|\oc_{\p^2}(\mu_1 - 1)-\sum_{i=2}^r \mu_i 
\pi(q_i)|$, whose virtual dimension (i.e. the dimension evaluated as if the points impose 
independent conditions) is given by the following formula:
\[
\binom{\mu_1 + 1}{2}-\sum_{i=2}^r\binom{\mu_i + 1}{2} - 1.
\]
By means of the two conditions $\langle C,C\rangle=-3$ and $\langle C,K\rangle=0$, we know 
that $\sum_{i=2}^r \mu_i^2 = \mu_1^2  - 2\mu_1 + 2$ and $\sum_{i=2}^r\mu_i = 3\mu_1 -2$ . 
In this way the preceding expression becomes equal to $-1$. 
We recall that Harbourne-Hirschowitz conjecture (see~\cite{har,hi}) implies that either the 
system $\ls_2$ 
is empty (and then also the set $\ell(2\mu_1-1;\mu_1,\ldots,\mu_r)$ is), or 
it is non-reduced. In this case, since $\Pi$ is generic then 
either $C$ is non-reduced or it contains a plane curve of degree $\geq 2$ through $q_1$
(which is projected to a multiple line).
\end{proof}

We end the section with two examples.

\begin{example}
A curve $C\in\ell(13;6,4^2,3,1^9)$ satisfies 
equalities~\eqref{-1} and $\delta > 2\mu_1$. Since $\delta-\sum_{i=1}^4\mu_i= - 4$, 
applying a transformation based on the first $4$ points we get 
a curve $C'\in\ell(7;2,-1,1^9) = \ell(7;2,1^9) + e_4$ 
(the class of a line contained in the exceptional divisor $E_4$). 
Therefore, since the transformation restricted to $E_4$ is a planar Cremona transformation, 
$C$ is reducible and contains the planar component $C_1\in\ell(2;1^3)$ (which goes to $e_4$). \\
In fact if we project $C$ from $p_1$ on a plane we get a curve in $\ls_2(7;4^2,3,1^9)$ which 
contains the line through the first two points. This is the image of the conic $C_1$ via the 
projection.\\
In particular, every curve $C\in\ell(13;6,4^2,3,1^9)$ is reducible and hence it is
not a numerically elementary $(-1)$-curve.
\end{example}

\begin{example}
Let us consider a curve $C\in\ell(7;4,1^{10})$. Following the proof of Proposition~\ref{proj} 
let us project $C$ from $q_1$ to the curve $\Gamma$ contained in a generic plane. 
The curve $C$ must contain the line $\langle q_1,q_i\rangle$, for some $i\in\{2,\ldots,11\}$,  
since otherwise $\Gamma$ would belong to the planar system 
$\ls_2(3;1^{10})$ which is empty. Therefore, as before $C$ must be reducible, which means
that $\ell(7;4,1^{10})$ contains no numerically elementary $(-1)$-curves.\\
Suppose now that the $10$ points lie on the cone $V$ over a plane cubic with vertex $q_1$. 
Blowing-up $V$ along $q_1$  we obtain a ruled elliptic surface $S$.
Let $H\in\pic(S)$ be the pull-back of an hyperplane section of the cone and let $F$ be the 
class of a fiber. The exceptional divisor $E$ is numerically equivalent to $H-3F$ and in 
particular $E^2=-3$. The linear system $|H+4F|$ is non-special, hence by Riemann-Roch it 
has dimension 10. Observe that the system can not be composed with a pencil since otherwise 
$H+4F\cong n(aH+bF)$, which can not happen if $n>1$. Since the system has no fixed components, 
then by Bertini's Theorem its general element is an irreducible elliptic curve. This implies that trough 10 
points in general position of $S$ there exists an irreducible curve of $|H+4F|$. The blow-down of $E$ is 
given by the map $\phi_{|H|}: S\rt V$ so, since $(H+4F) E=4$ and $(H+4F)H=7$, an element 
$\Gamma\in |H+4F|$ is sent to an irreducible elliptic curve of $\ell(7;4,1^{10})$. 
This  somewhat 
strange phenomenon is due to the fact that the points  are no longer in general position. 
This also implies that such a curve can not be obtained as a specialization of a general one.
\end{example}

\bibliographystyle{plain}

\end{document}